\documentclass[12pt]{article}
\usepackage{amsmath,amscd,amsbsy,amssymb,latexsym,url,bm,amsthm}
\usepackage[vlined,boxed,commentsnumbered,linesnumbered,ruled]{algorithm2e}
\usepackage{epsfig,graphicx,subfigure}
\usepackage{enumitem,balance,mathtools}
\usepackage{wrapfig}
\usepackage{mathrsfs, euscript}
\usepackage{amsmath, amssymb, amsthm}
\usepackage{epstopdf}
\usepackage{float}
\usepackage{comment}
\usepackage{abstract}
\usepackage{newclude}
\author{Jiyou Li, Yanghongbo Zhou}
\title{A Bijection between Necklaces and Restricted Multisets}
\date{}
\begin{document}
	\vspace{1.5cm}
	\newpage
	\renewcommand{\contentsname}{Content} 
	\newpage
\newtheorem{theorem}{Theorem}
\newtheorem{lemma}[theorem]{Lemma}
\newtheorem{proposition}[theorem]{Proposition}
\newtheorem{corollary}[theorem]{Corollary}
\newtheorem{definition}{Definition}
\newtheorem{remark}{Remark}[section]
\newtheorem{example}[theorem]{Example}
\newtheorem{exercise}{Exercise}
\newenvironment{solution}{\begin{proof}[Solution]}{\end{proof}}
	\setlength{\lineskiplimit}{2.625bp}
	\setlength{\lineskip}{2.625bp}
	\numberwithin{equation}{section}
	\newenvironment{partlist}[1][]
	{\begin{enumerate}[itemsep=0pt, label=(\arabic*), wide, labelindent=\parindent, listparindent=\parindent, #1]}
		{\end{enumerate}}
	\setcounter{page}{1}
	\maketitle
\begin{abstract}
We present a proof of Swee Hong Chan's conjecture establishing a bijection between 
the set of necklaces of length $n$ with at most $q$ colors, and the set of periodic functions $f: \mathbb{Z}_{n}\to {0, 1, ..., q-1}$ whose weighted sum is divisible by $n$, where $q$ and $n$ are coprime positive integers.

\end{abstract}

	\section{Introduction}

 The Subset Sum Problem (SSP) is a fundamental problem in theoretical computer science and combinatorial optimization. The standard SSP asks for the existence of a subset \( S \) in a finite set \( D \) of integers such that the sum of elements in \( S \) equals a given target. As one of the most well-known NP-complete problems, SSP has wide applications in complexity theory and cryptography.

While classically defined over the integers, its generalization to finite abelian groups reveals deeper structural insights and broader computational implications. This framework preserves core computational characteristics while highlighting algebraic structure, where SSP over prime-order groups remains NP-hard under RP-reduction as integer SSP can be embedded into these groups \cite{cheng2011countingvaluesetsalgorithm}.

Let \( G \) be an abelian group and \( D \subset G \) a finite subset of \( n \) elements. For an element \( b \in G \), let \( N_D(b) \) denote the number of subsets \( S \subseteq D \) satisfying \(\sum_{a \in S} a = b\).  The decision version is to determine existence of a non-empty subset summing to \( b \) (i.e., whether \( N_D(b) > 0 \)), while the counting version computes \( N_D(b) \) explicitly. For convenience, in the counting version the empty set is regarded as a set summing to 0. 

A particularly well-studied case arises when \( G = \mathbb{Z}/n\mathbb{Z} \) and \( D = G \),
 due to its connections to combinatorics and number theory.
  In this case the cardinality  $N_{\mathbb{Z}/n\mathbb{Z}}(0)$
admits the following elegant closed-form expressions:
 \[
N_{\mathbb{Z}/n\mathbb{Z}}(0) = \frac{1}{n} \sum_{\substack{d \mid n \\ d \text{ odd}}} \varphi(d) 2^{n/d},
\]
where $\varphi(d)$ is the Euler totient function. 
This formula also counts the number of inequivalent binary necklaces under the action of the cyclic group of order $n$. It can be derived from Burnside's lemma on orbits enumeration and reflects interesting combinatorial symmetry. Stanley \cite{MR4621625} constructed an explicit bijection between binary necklaces of length $n$ and zero-sum subsets of \(\mathbb{Z}/n\mathbb{Z} \) for prime $n$ and asked a bijective proof for non-prime $n$. Kaseory later extending the equivalence of cardinalities to arbitrary odd $n$ via private communication \cite{MR4621625}. For \( b \neq 0 \), \( N_{\mathbb{Z}/n\mathbb{Z}}(b) \) lacks a simple closed form and requires more intricate combinatorial techniques. See \cite{MR2457537,MR4621625}.

The correspondence deepens when considering $q$-colored necklaces. Chan's groundbreaking work \cite{MR3934368} established a polynomial-based bijection between:
\begin{itemize}
    \item $\mathcal{N}(n, q)$: $q$-colored necklaces with $n$ beads (equivalent up to cyclic rotations)
    \item $\mathcal{F}(n, q)$: Restricted multisets $S \subseteq \mathbb{Z}/n\mathbb{Z}$ satisfying $\sum_{i=1}^n i c_i \equiv 0 \pmod{n}$ with multiplicities $c_i < q$
\end{itemize}
through ingenious use of finite field decompositions and discrete logarithms, though limited to prime power $q$.

The recent work of Li and Yu \cite{MR4265614} established the enumerative equality for all coprime $(n,q)$, motivating our complete constructive proof of the bijection in this general case.  Separately, inspired by a bijection between complete linear systems and binary necklaces \cite{MR4159825}, Oh and Park \cite{MR4080626} discovered an alternative bijection (independent of Chan's algebraic approach) between binary necklaces with $n$ black and $k$ white beads and special $(n,k)$-error correcting codes valid when either $n$ or $k$ is prime. This further enriches the combinatorial landscape. 
	\section{Main Result}
\hspace{2em}
Let \( n \) and \( q \) be coprime positive integers. Suppose \( q \) has the prime factorization \( q = \prod_{i=1}^k q_i \), where each \( q_i = p_i^{t_i} \) is a prime power of a distinct prime \( p_i \). Our construction relies on the following foundational bijection established by Swee-Hong Chan:
\begin{theorem}\label{bijection_theorem}\cite{MR3934368}
Let \( q \) be a prime power and \( n \) a positive integer coprime to \( q \). Then there exists an explicit bijection
\[
\Phi_{n,q}: \mathcal{N}(n, q) \xrightarrow{\sim} \mathcal{F}(n, q).
\]
\end{theorem}

 The core innovation of our approach is to combine the bijections \( \Phi_{n,q_i} \) using the arithmetic identity:
\[
q - 1 = \sum_{i=1}^k (q_i - 1) \cdot \frac{q}{\prod_{j=1}^i q_j}.
\]
This identity allows us to synthesize global combinatorial data from the local bijections in Theorem~\ref{bijection_theorem}. To formalize this, we require the following decomposition:

 \begin{definition}[Cyclotomic cosets]
\label{def1}
For \( q = \prod_{i=1}^k q_i \) with \( q_i = p_i^{t_i} \) and distinct primes \( p_i \), define an equivalence relation on \( \mathbb{Z}/n\mathbb{Z} \) by the group action \( x \mapsto q_i \cdot x \). The equivalence classes (called \emph{cyclotomic cosets}) partition \( \mathbb{Z}/n\mathbb{Z} \) into disjoint sets:
\begin{align*}
&\mathscr{S}_{1,1},\ \mathscr{S}_{1,2},\ \dots,\ \mathscr{S}_{1,m_1}, \\
&\mathscr{S}_{2,1},\ \mathscr{S}_{2,2},\ \dots,\ \mathscr{S}_{2,m_2}, \\
&\qquad\vdots \\
&\mathscr{S}_{k,1},\ \mathscr{S}_{k,2},\ \dots,\ \mathscr{S}_{k,m_k}.
\end{align*}
Each coset has the form
\[
\mathscr{S}_{i,j} = \left\{ s_{i,j} \cdot q_i^\ell \mod n  \mid  \ell = 0,1,\dots,\lambda_{i,j}-1 \right\},
\]
where \( s_{i,j} \) is a minimal representative, and \( \lambda_{i,j} \) is the minimal positive integer satisfying 
\[
q_i^{\lambda_{i,j}} \cdot s_{i,j} \equiv s_{i,j} \pmod{n}.
\]
Here \( m_i \) denotes the number of distinct cosets for fixed \( i \).
\end{definition}

%

We view the set of necklaces \(\mathcal{N}(n, q)\) from the following algebraic perspective. For each prime power factor \(q_i\), define the quotient ring
\[
Q_i := \mathbb{F}_{q_i}[X]/(X^n - 1).
\]
By the Chinese Remainder Theorem for polynomial rings, there exists a ring isomorphism:
\begin{corollary} \label{Chinese_Remainder}\cite{MR2286236}
Let \(q = \prod_{i=1}^k q_i\) where each \(q_i = p_i^{t_i}\) is a prime power of distinct primes. Then there is a ring isomorphism:
\[
\mathbb{Z}/q\mathbb{Z}[X]/(X^n - 1) \cong \prod_{i=1}^k \mathbb{F}_{q_i}[X]/(X^n - 1).
\]
\end{corollary}

Via this isomorphism, any element \(\alpha \in \mathbb{Z}/q\mathbb{Z}[X]/(X^n - 1)\) corresponds to a \(k\)-tuple \((\alpha_1, \dots, \alpha_k)\) with \(\alpha_i \in Q_i\). Consequently, the necklace set admits the representation:
\[
\mathcal{N}(n, q) = \left\{ \left\langle \alpha \right\rangle \Bigm| \alpha \in \mathbb{Z}/q\mathbb{Z}[X]/(X^n - 1) \right\}
\]
where \(\langle \beta \rangle=\{X^j \cdot \beta\}_{j=0}^{n-1}\) denotes the equivalence class under cyclic shifts. In terms of the decomposition:
\[
\mathcal{N}(n, q) \cong \prod_{i=1}^k \mathcal{N}(n, q_i),
\]
with each \(\mathcal{N}(n, q_i)\) consisting of orbits
\[
\left\{ \left\langle \alpha_i \right\rangle \Bigm| \alpha_i \in Q_i \right\}.
\]

Since \(n\) and \(q_i\) are coprime, each field \(\mathbb{F}_{q_i}\) contains a primitive \(n\)-th root of unity \(\omega_i\) in its algebraic closure \(\overline{\mathbb{F}}_{q_i}\).

%
%

\begin{definition}\label{def2}
Let $q_{i}$ be a prime power as mentioned before. Let $P_{i,1},\cdots,P_{i,m_{i}}$ be the irreducible factors of $X^{n}-1$ over the field $\mathbb{F}_{q_{i}}$. That is, for any $j_{i}\in \left\{1,...,m_{i}\right\}$,
$$P_{i,j}:=\prod_{k\in \mathscr{S}_{i,j}}(X-\omega_{i}^{k}).$$
\end{definition}
We refer to \cite{MR2008834} for the proofs of the following properties of $\mathcal{Q}_{i}/P_{i,j}\mathcal{Q}_{i}$ and $G_{i,j}$.\\
\begin{lemma}\cite{MR2008834}
 Let $q_{i}$ be a prime power, and let $n$ be a positive integer coprime to $q_{i}$. For any $j_i\in \{1,\cdots, m_i\}$,
 \begin{itemize}
     \item $\mathcal{Q}_i/P_{i,j}\mathcal{Q}_i$ is a finite field of order $q_{i}^{\lambda_{i,j}}$.
     \item $G_{i,j}$ is a cyclic group of order $q_{i}^{\lambda_{i,j}}-1$ order under multiplication.
     \item $X_{i,j}$ is an element of $G_{i,j}$ with multiplicative order $\frac{n}{\text{gcd}(n,s_{i,j})}$.
 \end{itemize}
\end{lemma}
The following lemma is based on the Chinese remainder theorem. 
\begin{lemma}\cite{MR3934368}
 Let $n$ and $d_1,\cdots,d_{k}$ be positive integers. Then there exists a group automorphism $\phi:\prod^{k}_{i=1}\mathbb{Z}_{\frac{n}{\text{gcd}(n,d_{i})}}\to\prod^{k}_{i=1}\mathbb{Z}_{\frac{n}{\text{gcd}(n,d_{i})}}$ such that
 \begin{equation*}
   d_1h_1+\cdots+d_kh_k=\text{gcd}(n,d_1,\cdots,d_{k}) \mod n
 \end{equation*}
 where $h_i$ is the $i$-th coordinate of $\phi(1,\cdots,1)$.
\end{lemma}

Though Jiyou Li and Xiang Yu have gave them skillful proof on the sizes of the two objects, we will present another method.
\subsection{Proof of Cardinalities}
\hspace{2em}For any function $f\in\mathcal{F}(n,q)$, define the level set $L_{q-1}(f)$ be the set $\{z\in\mathbb{Z}_{n}|f(z)=q-1\}$.
\begin{definition}\label{def4}
Let $q$ be a number coprime to $n$. For any $I_{i}\subseteq \left\{1,\cdots,m_{i}\right\}$, the sets $\mathcal{N}(n,q)_{I_{1},\cdots,I_{k}}$ and $\mathcal{F}(n,q)_{I_{1},\cdots,I_{k}}$ are defined as follows.  
$$
\begin{aligned}
\mathcal{N}(n,q)_{I_{1},\cdots,I_{k}}&:=\left\{\prod^{k}_{i=1}\left\langle \alpha_i \right\rangle\in\mathcal{N}(n,q)|\alpha_{i}\not\equiv0\mod P_{i,j}\right\}\\
\mathcal{F}(n,q)_{I_{1},\cdots,I_{k}}&:=\left\{f\in\mathcal{F}(n,q)|L_{q-1}(f)\cap \mathscr{S}_{i,j}=\mathscr{S}_{i,j}\text{ iff }j\not\in I_{i}\right\}.
\end{aligned}
$$
\end{definition}
By definition, $\left\{\mathcal{N}(n,q)_{I_1,\cdots,I_k}\right\}_{I_i\subseteq\left\{1,\cdots,m_i\right\}}$ and $\left\{\mathcal{F}(n,q)_{I_1,\cdots,I_k}\right\}_{I_i\subseteq\left\{1,\cdots,m_i\right\}}$ form a partition of $\mathcal{N}(n,q)$ and $\mathcal{F}(n,q)$, respectively.\\
\begin{lemma}
Let $q$ be an integer, let $n$ be a positive integer coprime to $q$, and let $I_{i}\subseteq\left\{1,\cdots,m_i\right\}$. Then
$$|N(n,q)_{I_1,\cdots,I_k}|=\frac{gcd(n,(gcd(s_{i,j})_{j\in I_{i}})_{i=1}^{k})}{n}\prod_{i=1}^{k}\prod_{j\in I_i}(q_{i}^{\lambda_{i,j}}-1).$$
\end{lemma}
\begin{proof}
Recall the definition of $G_{i,j}$ from Definition \ref{def2} and the definition of $\alpha_{i,j}$ and $X_{i,j}$ from  previous definition. In particular, if $\alpha_{i}$ is an element of $\mathcal{Q}$ that is not divisible by $P_{i,j}$, then $\alpha_{i,j}$ is contained in $G_{i,j}$. Consider the map
$$
\begin{aligned}
\xi:\prod_{i=1}^{k}\left\{\alpha_{i}\in\mathcal{Q}|P_{i,j}\text{ divides }\alpha_{i}\text{ iff }i\not\in I\right\}
&\to
\prod_{i=1}^{k}\prod_{j\in I_{i}}G_{i,j}\\
(\alpha_{i})^{k}_{i=0} &\mapsto ((\alpha_{i,j})_{j\in I_i})^{k}_{i=0}
\end{aligned}
$$
The map $\xi$ is a bijection by Chinese Remainder Theorem in polynomial rings.\\
Denote by $C_{I_1,\cdots,I_m}$ the product of cyclic subgroups of $\prod_{i=1}^{k}\prod_{j\in I_i}G_{i,j}$ generated by $(X_{i,j})_{j\in I_i, 1\leq i\leq k}$. Note that $\mathcal{N}_{I_1,\cdots,I_k}$ is in bijection with cosets of $C_{I_1,\cdots,I_k}$ in $\prod_{i=1}^{k}\prod_{j\in I_i}G_{i,j}$ by the map $\xi_{i}$. Hence we have
\begin{equation}\label{eq1}
|\mathcal{N}_{I_1,\cdots,I_k}|=
|\prod^{k}_{i=1}\prod_{j\in I_i} G_{i,j}/C_{I_1,\cdots,I_k}|
=\frac{1}{C_{I_1,\cdots,I_k}}\prod^{k}_{i=1}\prod_{j\in I_i}|G_{i,j}|.
\end{equation}
On the other hand we also have
$$
\begin{aligned}
|G_{i,j}|&=q_{i}^{\lambda_{i,j}}-1\\
|C_{I_1,\cdots,I_k}|&=\min\left\{\tau>0|(X_{i,j})^{\tau}\text{ is the identity element of } G_{i,j}\text{ for all }j\in I_i, 1\leq i\leq k\right\}\\
&=\operatorname{lcm}(1,(\frac{n}{\gcd(n,s_{i,j})_{j\in I_i}}))\\
&=\frac{n}{\gcd(n,\gcd(n,s_{i,j})_{j\in I_i})}
\end{aligned}
$$
The conclusion of the lemma now follows from (\ref{eq1}).
\end{proof}

\begin{lemma}\label{size}
Let $q$ and $n$ be two comprime positive integers, and let $I_{i}\subseteq\{1,\cdots,m_{i}\}$. Then
\begin{equation}
|\mathcal{F}(n,q)_{I_1,\cdots,I_{k}}|=\frac{gcd(n,gcd(s_{i,j})_{j\in I_{i}})}{n}\prod_{i=1}^{k}\prod_{j\in I_i}(q_{i}^{\lambda_{i,j}}-1).
\end{equation}
\end{lemma}
\begin{proof}
Let $\varepsilon_{I_{1},\cdots,I_{k}}$ denote the set
\begin{equation}
\varepsilon_{I_1,\cdots,I_{k}}:=\{f:\mathbb{Z}_{n}\to\{0,1,\cdots,q-1\}|L_{q-1}(f)\cap \mathscr{S}_{i,j}=\mathscr{S}_{i,j}\text{ iff }j\not\in I_{i}\}.
\end{equation}
Let $\eta_{I_1,\cdots,I_k}:\varepsilon_{I_{1},\cdots,I_{k}}\to\prod_{i}\prod_{j\in I_{i}}$ be the map defined by
\begin{equation}
f\mapsto \left(\sum^{\lambda_{i,j}-1}_{u=0}f(q_{i}^{u}s_{i,j}) \mod q_{i}^{\lambda_{i,j}-1}\right)_{j\in I_{i},i\in\{1,\cdots,k\}}.
\end{equation}
The map is a surjective by the definition of $\varepsilon_{I_1,\cdots,I_{k}}$.
Let $f$ be any function in $\varepsilon_{I_{1},\cdots,I_{k}}$. For any $j\in I_{i}$, the corresponding sum $\sum^{\lambda_{i,j}-1}_{u=0}f(q_{i}^{u}s_{i,j})$ is strictly less than $q_{i}^{\lambda_{i,j}}-1$ since $L_{q-1}\cap \mathscr{S}_{i,j}\not= \mathscr{S}_{i,j}$. This imples that the $i$-th coordinate of $\eta_{I_1,\cdots,I_{k}}(f)$ determines $f(s_{i,j}),\cdots,f(s_{i,j}q_{i}^{\lambda_{i,j}-1})$ for any $j\in I_{i}$. Furthermore, we have $f(s_{i,j}),\cdots,f(s_{i,j}q_{i}^{\lambda_{i,j}-1}) = (q_{i}-1,\cdots,q_{i}-1)$ for any $i\not\in I$ by the definition of $\varepsilon_{i_1,\cdots,I_{k}}$. Therefore, we conclude that $\eta_{I_1,\cdots,I_{k}}$ is a surjective map.
Let $\zeta_{I_1,\cdots,I_{k}}$ be the map defined by
\begin{equation}
\begin{aligned}
\zeta_{I_{1},\cdots,I_{k}}:&\prod_{i=1}^{k}\prod^{j\in I_{i}}\mathbb{Z}_{q_{i}^{\lambda_{i,j}-1}}&\to&\mathbb{Z}_{n}\\
&((z_{j})_{j\in I_{i}})_{i\in \{1,\cdots,k\}}&\mapsto&\sum_{i=1}^{k}\frac{q}{\prod_{u\leq i}q_{u}}\sum_{j\in I_i}s_{i,j}z_{i,j} \mod n
\end{aligned}
\end{equation}
The map $\zeta_{I_{1},\cdots,I_{k}}$ is a well defined group homomorphism as $n$ divides $s_{i,j}(q^{\lambda_{i,j}}-1)$ by previous definition. Further, by the definition of $\text{gcd}$, the image of $\zeta_{I_1,\cdots,I_{k}}$ is $\text{gcd}(n,\text{gcd}(s_{i,j})_{j\in I_{i}})\mathbb{Z}_{n}$.
Now note that, for any $f\in \varepsilon_{I_{1},\cdots,I_{k}}$,
\begin{equation}
\begin{aligned}
\sum_{z\in \mathbb{Z}_{n}}zf(z)&=\sum_{i=1}^{k}\frac{q}{\prod_{u\leq i}q_{u}}(\sum_{j_i\in I_i}\sum^{l_{i,j_i}-1}_{u=0}q_{i}^{u}s_{i,j_i}f(q^{u}s_{i,j_i})+\sum_{j_i\not\in I_{i}}(q_{i}^{l_{i,j_i}}-1)s_{i,j_i})\\
&=\sum_{i=1}^{k}\frac{q}{\prod_{u\leq i}q_{u}}(\sum_{j_i\in I_i}s_{i,j_i}\sum^{l_{i,j_i}-1}_{u=0}q_{i}^{u}f(q_{i}^{u}s_{i,j_i}))\mod n\\
&=\zeta_{I_1,\cdots,I_k}(\eta_{I_1,\cdots,I_k}(f)).
\end{aligned}
\end{equation}
Since $\eta_{I_1,\cdots,I_{k}}$ is a bijection, it then follows from the definition of $\mathcal{F}_{I_1,\cdots,I_{k}}$ that the kernel of $\zeta_{I_1,\cdots,I_{k}}$ is equal to $\eta_{I_1,\cdots,I_{k}}(\mathcal{F}(n,q)_{I_1,\cdots,I_{k}})$.
Combining all those observations, we conclude that
\begin{equation}
\begin{aligned}
|\mathcal{F}(n,q)_{I_1,\cdots,I_k}|&=|\eta_{I_1,\cdots,I_{k}}(\mathcal{F}(n,q)_{I_1,\cdots,I_k})| = |ker(\zeta_{I_1,\cdots,I_{k}})|\\
&=\frac{|\prod_{i=1}^{k}\prod_{j\in I_i}\mathbb{Z}_{q_{i}^{\lambda_{i,j}}-1}|}{|\text{gcd}(n,\text{gcd}(s_{i,j})_{j\in I_i, i\in\{1,\cdots,k\}})|}\\
&=\frac{\text{gcd}(n,\text{gcd}(s_{i,j})_{j\in I_i, i\in\{1,\cdots,k\}})}{n}\prod^{k}_{i=1}\prod_{j\in I_i}(q_{i}^{\lambda_{i,j}}-1)
\end{aligned}
\end{equation}
as desired.
\end{proof}
Combine the above two lemmas, we can get the equivalence of the sizes.

\subsection{Proof of Bijection}
In this section, we present a proof of the following Theorem:
\begin{theorem}
Let $n,q$ be two coprime positive integers. Then there is a bijection between $\mathcal{N}(n,q)$ and $\mathcal{F}(n,q)$.
\end{theorem}
Throughout this section, $q_{i}=p_{i}^{t_i}$ are prime powers for different primes and $n$ is a positive integers that is coprime to any of $q_{i}$.\\
Let $\mathcal{Q}_{i}$ be as defined before, and let $\mathcal{E}$ be the set of all functions from $\mathbb{Z}_{n}$ to $\left\{0,1,\cdots,q-1\right\}$. Suppose that there exists a map $\psi:\prod\mathcal{Q}_{i}\to\mathcal{E}$ that satisfies the following conditions:
\begin{itemize}
    \item The map $\psi$ is a bijection from $\mathcal{Q}$ to $\mathcal{E}$,
    \item For any $\alpha_{i} \in \mathcal{Q}_{i}$ there exists a unique $\beta\in\left\{\alpha_{i},X\alpha,\cdots,X^{n-1}\alpha\right\}$.
\end{itemize}
We could then define the map $\hat{\psi}:\mathcal{N}(n,q)\to\mathcal{F}(n,q)$ by 
$$
\left\{\alpha, X\alpha, \cdots, X^{n-1}\alpha\right\} \mapsto \psi(\beta).
$$
It would follow that $\hat{\psi}$ is a bijection between $\mathcal{N}(n,q)$ and $\mathcal{F}(n,q)$, which would prove our main result. In this section, we will construct a map $\psi:\mathcal{Q}\to\mathcal{E}$ that satisfies the two conditions.\\
Recall the definition of $m_{i}, s_{i,j}$ and $\lambda_{i,j}$ from Definition \ref{def1}, the definition of $G_{i,j}$ from Definition \ref{def2} and the definition of $X_{i,j}$ from previous definition.\\
Let $j \in \left\{1, \cdots, m_{i}\right\}$. Since $G_{i,j}$ is a cyclic group of order $q_{i}^{\lambda_{i,j}}$ and $X_{i,j}$ is an element of $G_{i,j}$ with order $\frac{n}{\gcd(n,s_{i,j})}$, the group $G_{i,j}$ contains a group generator such that $X_{i,j}$ is $\frac{q_{i}^{\lambda_{i,j}}\gcd(n,s_{i,j})}{n}$-th power of this generator.\\
\begin{definition}\label{def5}
For any $j\in \left\{1,\cdots,m_{i}\right\}$, let $g_{i,j}$ be a group generator of $G_{i,j}$ such that $X_{i,j}$ is the $\frac{q_{i}^{l_{i,j}}\gcd(n,s_{i,j})}{n}$-th power of $g_{i,j}$.
\end{definition}
Recall the definition of $P_{i,j}$ from Definition \ref{def2} and the definition of $\alpha_{i,j}$ from previous definition.
\begin{definition}\label{def6}
Let $j \in \left\{1,\cdots,m_i\right\}$, and let $\alpha$ be an element  of $\mathcal{Q}$ not divisible by $P_{i,j}$. The discrete logarithm $\log_{g_{i,j}}(\alpha)$ is the smallest non-negative integer $k$ such that $\alpha_{i,j}=g_{i,j}^k$. 
\end{definition}
By the lemma of cyclic group, the integer $\log_{g_{i,j}}(\alpha)$ is contained in $\left\{0,\cdots,q_{i}^{\lambda_{i,j}}-2\right\}$.
\begin{definition}\label{def7}
Let $j\in\left\{1,\cdots,m_i\right\}$, and let $\alpha$ be an element of $\mathcal{Q}$ not divisible by $P_{i,j}$. We denote by $a_{i,j}(\alpha)$ and $b_{i,j}(\alpha)$ the quotient and the remainder of the division of $\log_{g_{i,j}}(\alpha)$ by $\frac{q_{i}^{\lambda_{i,j}}\gcd(n,s_{i,j})}{n}$, respectively.
\end{definition}
\begin{lemma}
Let $q_i$ be a prime power, let $n$ be a positive integer coprime to $q_i$ and let $j \in \left\{1,\cdots,m_i\right\}$. Then
\begin{enumerate}
    \item\label{prt1} $a_{i,j}(X)=1$ and $b_{i,j}(X)=0$,\\
    \item\label{prt2} For any $k \geq 0$ and any $\alpha \in \mathcal{Q}$,
          \begin{equation*}
          \begin{aligned}
            a_{i,j}(X^{k}\alpha) &= k+a_{i,j}(\alpha) (\mod \frac{n}{gcd(n,s_{i,j})})\\
            b_{i,j}(X^{k}\alpha) &= b_{i,j}(\alpha).
            \end{aligned}
          \end{equation*}    
\end{enumerate}
\end{lemma}
\begin{proof}
Part \ref{prt1} follows directly from Definition \ref{def5} and Definition \ref{def7}.\\
By Definition \ref{def6}, we have for any non-negative integer $k$ and any $\alpha \in \mathcal{Q}$ that
  \begin{equation*}
    \begin{aligned}
      \log_{g_{i,j}}(X^{k}\alpha)&=\log_{g_{i,j}}(\alpha)+k\log_{g_{i,j}}(X) (\mod q_{i}^{l_{i,j}}-1)\\
      &=(k+a_{i,j}(\alpha))\frac{q_{i}^{l_{i,j}}\gcd(n,s_{i,j})}{n}+b_{i,j}(\alpha) (\mod q_{i}^{\lambda_{i,j}}-1).
    \end{aligned}
  \end{equation*}
Part \ref{prt2} now follows from Definition \ref{def7}.
\end{proof}

\begin{definition}
  Let $I_{i}$ be a subset of $\left\{1,\cdots,m_i\right\}$. Let $\phi_{i,I_{i}}$ be a group automorphism of $\prod_{j\in I_{i}}\mathbb{Z}_{\frac{n}{\gcd(n,s_{i,j})}}$ that satisfies
  \begin{equation*}
    \sum_{i=1}^{k}\frac{q}{\prod_{u\leq i}q_{u}}\sum_{j\in I_{i}}s_{i,j}h_{i,j,I_i}=\gcd(n,\gcd(s_{i,j})_{j\in I_{i}, 1\leq i \leq k}) \mod n,
  \end{equation*}
  where $h_{i,j,I_i}$ is the $i,j$-th coordinate of $\phi_{i,I_i}(1,\cdots,1)$. The function $\phi_{i,I_i}$ exists for any $I_i\subseteq \left\{1,\cdots,m_i\right\}$ by a lemma.
\end{definition}
Recall that $L_{q-1}(f)=\left\{z\in \mathbb{Z}_{n}|f(z)=q-1\right\}$. For any $I_{i}\subseteq\left\{1,\cdots,m_i\right\}$, write
\begin{equation*}
  \begin{aligned}
    \mathcal{Q}_{i,I_i}&:=\left\{\alpha\in\mathcal{Q}|P_{i,j}\text{ divides } \alpha\text{ iff }j\not\in I_{i}\right\}\\
    \mathcal{E}_{i,I_i}&:=\left\{f\in\mathcal{E}|L_{q-1}(f)\cap S_{i,j}=S_{i,j}\text{ iff }j\not\in I_{i}\right\}.
  \end{aligned}
\end{equation*}
By definition $\prod(\mathcal{Q}_{i,I_i})_{I_i\subseteq [m_{i}]}$ and $\prod(\mathcal{E}_{i,I_i})_{I_i\subseteq [m_{i}]}$ form a partition of $\mathcal{Q}$ and $\mathcal{E}$, respectively.\\
Let $j\in I_{i}$, and let $\alpha$ be any element of $\mathcal{Q}_{i,I_{i}}$. We denote by $\phi_{i,j,I_i}(\alpha)$ the $j$-th coordinate of $\phi_{i,I_{i}}(a_{i,j}(\alpha)_{j\in I_{i}})$, which corresponds to a non-negative integer strictly less than $\frac{n}{\gcd(n,s_{i,I_i})}$.\\
Since $b_{i,j}(\alpha)$ is a non-negative integer strictly less than $\frac{(q_{i}^{\lambda_{i,j}}-1)\gcd(n,s_{i,j})}{n}$ and $\phi_{i,j,I_i}(\alpha)$ is a non-negative integer strictly less than $\frac{n}{\gcd(n,s_{i,j})}$, we have
\begin{equation}\label{eq4}
  0 \leq b_{i,j}(\alpha)\frac{n}{\gcd(n,s_{i,j})}+\phi_{i,j,I_i}(\alpha)<q_{i}^{\lambda_{i,j}}-1. 
\end{equation}
We denote by $c_{i,j,0}(\alpha), \cdots, c_{i,j,\lambda_{i,j}-1}(\alpha)\in\left\{0,\cdots,q_{i}-1\right\}$ the unique integers that satisfy
\begin{equation}\label{eq5}
  \sum^{l_{i,j}-1}_{u=0}c_{i,j,u}(\alpha)q_{i}^{u}=b_{i,j}(\alpha)\frac{n}{\gcd(n,s_{i,j})}+\phi_{i,j,I_i}(\alpha).
\end{equation}
By \ref{eq4}, the sequence of integers $(c_{i,j,0}(\alpha), \cdots, c_{i,j,\lambda_{i,j}-1}(\alpha))$ is well defined and is not equal to $(q_{i}-1,\cdots,q_{i}-1)$.\\
Let $f_{\alpha_{i}}:\mathbb{Z}_{n} \to \{0, 1, \cdots, q_{i}-1\}$ be given by 
\begin{equation}\label{eq6}
    f_{\alpha_{i}}(q_{i}^{u}s_{i,j}):=\left\{
    \begin{aligned}
      q_{i}-1 & \text{ if }j\not\in I_{i};\\
      c_{i,j,u} & \text{ if }j\in I_{i}.
    \end{aligned}
    \right.
\end{equation}
$f_{\alpha}$ is then given by $f_{\alpha}(v)=\sum_{i=1}^{k}\frac{q}{\prod_{u\leq i}q_{u}}f_{\alpha_{i}}(v)$.
\begin{lemma}
  Let $q, n$ be two coprime numbers, and let $I\subseteq [k]$ and $J_i \subseteq \{1,\cdots, m_{i}\}$. Then for any $\alpha \in \prod_{i\in I, j\in J_{i}}{Q_{I_{i}}}$
  \begin{itemize}
      \item $\sum_{z\in \mathbb{Z}_{n}}zf_{\alpha}(z)=\sum_{i\in I}\frac{q}{\prod_{u\leq i}q_{u}}\sum_{j\in J_{i}}s_{i, j}\phi_{i,j,J_{i}} \mod n$;
      \item There exists unique $\beta \in \{\alpha, X\alpha, \cdots, X^{n-1}\alpha\}$ such that $\psi_{I, (J_{i})_{i\in I}}(\beta)$ is contained in $\mathcal{F}$.
  \end{itemize}
\end{lemma}
\begin{proof}
  \begin{itemize}
      \item  The first part is similar to the proof of lemma \ref{size}.
      \begin{equation*}
      \begin{aligned}
      \sum_{z\in \mathbb{Z}_{n}}z f_{\alpha}(z) &=\sum_{i\in I}\frac{q}{\prod_{u\leq i}q_{u}}\sum_{j\in J_{i}}\sum_{u=0}^{\lambda_{i, j}-1}q_{i}^{u}s_{i, j}c_{i, j}(\alpha)+\sum_{i\in I}\frac{q}{\prod_{u\leq i}q_{u}}\sum_{j\not\in J_{i}}\sum_{u=0}^{\lambda_{i, j}-1}q_{i}^{u}s_{i, j}(q_{i}-1) \mod n\\
      &=\sum_{i\in I}\frac{q}{\prod_{u\leq i}q_{u}}\sum_{j\in J_{i}}\sum_{u=0}^{\lambda_{i, j}-1}q_{i}^{u}s_{i, j}c_{i, j}(\alpha) \mod n\\
      &=\sum_{i\in I}\frac{q}{\prod_{u\leq i}q_{u}}\sum_{j\in J_{i}}(b_{i,j}(\alpha)\frac{n}{\gcd(n,s_{i,j})}+\phi_{i,j,J_i}(\alpha)) \mod n\\
      &=\sum_{i\in I}\frac{q}{\prod_{u\leq i}q_{u}}\sum_{j\in J_{i}}s_{i, j}\phi_{i,j,J_{i}} \mod n.
    \end{aligned} 
  \end{equation*}
   This proves the first part.
   \item We have $|\{\alpha, X\alpha, \cdots, X^{n-1}\alpha\}| = \text{lcm}(1, (\text{order of }X\text{ in }G_{i, j})) = \frac{n}{\text{gcd}(n, (s_{i, j})_{i\in I, j\in J_{i}})}$.
    We then have, for $k\geq 0$,
    \begin{equation*}
      \begin{aligned}
        \sum_{z\in \mathbb{Z}_{n}}zf_{X^{k}\alpha}(z)&=\sum_{i\in I}\frac{q}{\prod_{u\leq i}q_{u}}\sum_{j\in J_{i}}s_{i, j}\phi_{i, j, J_{i}}(X^{k}\alpha) \mod n\\
        &=\sum_{i\in I}\frac{q}{\prod_{u\leq i}q_{u}}\sum_{j\in J_{i}}s_{i, j}(kh_{i, j, J_{i}}+\phi_{i, j, J_{i}}(\alpha)) \mod n\\
        &=k\text{ gcd}(n, (s_{i, j})_{i\in I, j\in J_{i}}) + \sum_{i\in I}\frac{q}{\prod_{u\leq i}q_{u}}\sum_{j\in J_{i}}s_{i, j}\phi_{i, j, J_{i}}(\alpha) \mod n.
      \end{aligned}
    \end{equation*}
    By definition, $\sum_{i\in I}\frac{q}{\prod_{u\leq i}q_{u}}\sum_{j\in J_{i}}s_{i, j}\phi_{i, j, J_{i}}(\alpha)$ is a multiple of $\text{ gcd}(n, (s_{i, j})_{i\in I, j\in J_{i}})$, which completes the proof.
  \end{itemize}
\end{proof}

        \section{Example}
Let $q=10$ and $n=3$. We will show in this section how to compute the coresponding function of necklace $1+x+x^2$.
Since $q=2*5$, we suppose that $q_1=5$, $q_2=2$. Then the cyclotomic cosets should be:
\begin{equation*}
  \begin{aligned}
   s_{1,1}=0, & \mathscr{S}_{1,1}=\{0\},\\
   s_{1,2}=1, & \mathscr{S}_{1,2}=\{1, 2\}\\
   s_{2,1}=0, & \mathscr{S}_{2,1}=\{0\}\\
   s_{2,2}=1, & \mathscr{S}_{2,2}=\{1, 2\}.
  \end{aligned}
\end{equation*}
The irreducible polynomials $P_{i, j_{i}}$ are:
\begin{equation*}
  \begin{aligned}
    P_{1,1}=1+x, &\quad P_{1,2} = 1+x+x^2\\
    P_{2,1}=1+x, &\quad P_{2,2} = 1+x+x^2.
  \end{aligned}
\end{equation*}
We make the following choices of $g_{i, j_{i}}$ that satisfy the condition in Definition \ref{def5}.
\begin{equation*}
  \begin{aligned}
    g_{1,1}=1 \mod 1+x, &\quad g_{1,2} = x \mod 1+x+x^2\\
    g_{2,1}=1 \mod 1+x, &\quad g_{2,2} = x \mod 1+x+x^2.
  \end{aligned}
\end{equation*}
By Chinese Remainder Theorem,
\begin{equation*}
  \begin{aligned}
    Z_{10}(x)/(x^3-1)&\to \mathbb{F}_{2}(x)/(x^3-1)\times\mathbb{F}_{5}(x)/(x^3-1)&\to &P_{1,1}\times P_{1,2}\times P_{2,1}\times P_{2,2}\\
    \alpha=1+x+x^2 &\mapsto (1+x+x^2, 1+x+x^2)&\mapsto&(1,0,1,0).
  \end{aligned}
\end{equation*}
From Definition \ref{def7},
\begin{equation*}
  \begin{aligned}
    \log_{g_{1,1}}(\alpha)=0,& \log_{g_{2,1}}(\alpha)=0 
  \end{aligned}
\end{equation*}
For $\alpha$, $I_{1}=\{1\}, I_{2}=\{1\}$. Thus we want $s_{1,1}h_{1,1,1}+a_{2,1}h_{2,1,1}=3\mod 3$. Let $h_{1,1,1}=1$ and $h_{2,1,1}=1$, we get $0\cdot1+0\cdot1=3\mod3$.\\
Thus
\begin{equation*}
  \begin{aligned}
   b_{1,1}(\alpha)+\phi_{1,1\{1\}}=0+0=0\cdot5^{0}\\
   b_{2,1}(\alpha)+\phi_{2,1\{1\}}=0+0=0\cdot2^{0}.
  \end{aligned}
\end{equation*}
Both send $\{1+x+x^2\}$ to the function $0*0+0*1+0*2$. Thus the corresponding function of the necklace $\{1+x+x^2\}$ is $0*0+0*1+0*2$.
\bibliographystyle{plain}

\end{document}